   \input amstex
\documentstyle{amsppt}
\magnification=\magstep1
\parindent=1em
\CenteredTagsOnSplits
\NoBlackBoxes
\nopagenumbers
\NoRunningHeads
\pageno=1
\footline={\hss\tenrm\folio\hss}

        \topmatter
        \title {
     On linear operators with ${\ssize\bold p}$-nuclear adjoints
               }
        \endtitle
          \author { O.I.Reinov${{ }^\dag}$}  \endauthor
\address\newline
Oleg I. Reinov \newline
Department of Mathematics\newline
St Petersburg University\newline
St Peterhof, Bibliotech pl 2\newline
198904  St Petersburg, Russia
\endaddress

\email
orein\@orein.usr.pu.ru
\endemail

\thanks
${{ }^\dag}$This work was done with partial support of the Ministry of the
general and professional education of Russia (Grant 97-0-1.7-36) and
FCP ``Integracija", reg. No. 326.53.
\endthanks

\abstract
If $ p\in [1,+\infty]$ and $ T$ is a linear operator with $ p$-nuclear
adjoint from a Banach space $ X$ to a Banach space $ Y$ then if one
of the spaces $ X^*$ or $ Y^{***}$ has the approximation property, then
$ T$ belongs to the ideal $N^p$ of operators which can be factored
through diagonal oparators $l_{p'}\to l_1.$
On the other hand,
there is a Banach space $ W$ such that $ W^{**}$ has
a basis and such that for each $ p\in [1,+\infty], p\neq 2,$
there exists an operator $ T: W^{**}\to W$ with
$ p$-nuclear adjoint that is not in the ideal $N^p,$
as an operator from $ W^{**}$ to $ W.$
\endabstract
    \endtopmatter

\document
\baselineskip=18pt

\footnote""{${ }^\ddag$
}
\footnote""{${ }$
}

\def\al{\alpha}
\def\ot{\otimes}
\def\wh{\widehat}
\def\ffi{\varphi}
\def\wt{\widetilde}
\def\small{\smallpagebreak}
\def\sbs{\subset}
\def\({\left(}
\def\){\right)}
\def\({\left(}
\def\){\right)}
\def\[{\left[}
\def\]{\right]}
\def\tr{\operatorname{trace}\,}

In a series of previous notes the author considered the question
on the existence or nonexistence
of non-$p$-nuclear operators with $ p$-nuclear second
adjoints in concrete Banach spaces. For example, in [3] such
"bad" operators were constructed for every $ p\in [1,2)$
in spaces with Schauder bases.

Here we treat the analogous case for the operators
{\it whose first adjoints are $ p$-nuclear}. Already in the paper [2, \S4]
the corresponding assertions were formulated (without proofs)
in which the spaces under consideration, however, did not possesses
bases, --- in the best case only one of the spaces, in which
the operators acted, had the approximation property, but yet
not bounded AP. In this note, proving partly
those assertions from [2], we show that the similar examples 
can be found
in spaces with bases (and even in better spaces). Simultaneously,
some sufficient conditions are given for the positive answer to the
corresponding question. These questons are connected with the $p$-nuclearity
of the adjoint operator (moreover, those conditions are shown to be
close to the necessary ones).
\small

All the spaces under consideration are Banach spaces.
We usually denote the elements of the spaces by the corresponding small
letters: $x\in X,\,y\in Y,\dotsc,\, x'\in X^*,\,y''\in Y^{**},\dotso$.
For  $ p\in [1,+\infty],$
the conjugate exponent  $ p'$ is defined by the relation $ 1/p+1/p'=1.$
By $L\left(X,Y \right)$ it is denoted the space of all
(linear continuous) operators from $X$ to $Y$ with its standard norm.
We shall consider every Banach space also as a subspace in its
second dual, usually without introducing any additional notations
for the natural imbedding $ X\to X^{**}.$
However, in the case of a necessity,
we will denote by $ \pi_X$ this canonical imbedding.

Let us recall that {\it a Banach space has the {\rm(}Grothendieck{\rm)}
approximation
property}\ (the property AP), if the identity map on the space
can be approximated, in the topology of compact convergence,
by finite-dimensional operators.

For $p\in[1,\infty],$  an operator $T,$ acting from $X$ to $Y,$
is said to be
an {\it $N^p$-operator,} if it can be represented in the following form:
$$
  Tx= \sum_{k=1}^\infty \,<x,x'_k,>\,y_k \qquad\text{ for $x\in X$},
\tag1
$$
where the sequences
$\,\{x'_n\}^{\infty}_{n=1}\subset X^* \,$ and $\,\{y_n\}^{\infty}_{n=1}
\subset Y \,$ are such that the quantity
$$
  \al:=\, \left(\sum_{k=1}^\infty \,\|y_k\|^p\right)^{1/p}
    \sup \left\{ \left(\sum_{k=1}^\infty \,|<x'_k,x>|^{p'}\right)^{1/p'}:
     \  x\in X,\, \|x\|\le 1\right\}  \tag2
$$
is finite
(recall that $p'$ is the conjugate exponent for $p;$
in the case when one of the exponents $\,p,\,p'$ is equal to infinity,
one must understand the right part of the relation (2) properly).
The set of all $N^p$-operators from $X$ to $Y$ is denoted by
${N}\!{^p}\,(X,Y),$ and the lowest bound $\inf \al,$
where the {\it inf}\  is taken over all possible representations
in the form (1)
of the operator $T,$ --- by
$\nu^p\,(T).$ For every $p\ge 1$ the class ${N}\!{^p}$ of all
$N^p$-operators is a Banach operator ideal [1]; for fixed spaces
$X,Y$ \ \, ${N}\!{^p}\,(X,Y)$ is a Banach space with the norm
$\nu^p.$ Let us note that for the exponents
$p,q\in [1,+\infty],p<q,$ we always have:
${N}\!{^p} \subset {N}^q$ with the corresponding inequality for norms.
\small

We denote by $ X^*\wh\ot^p Y$
the tensor product associated with the space
${N}\!{^p}\,(X,Y)\,:$
this is the completion of the algebraic tensor product
$ X^*\ot Y$ via the norm $ \nu_0^p,$ which is defined for
$ z\in X^*\ot Y$ as the lowest bound, over all possible representations
$ z=\sum_{k=1}^N \,x'_k\ot\,y_k $ in $ X^*\ot Y,$ of the numbers
$$
   \left(\sum_{k=1}^N \,\|y_k\|^p\right)^{1/p}
    \sup \left\{ \left(\sum_{k=1}^N \,|<x'_k,x>|^{p'}\right)^{1/p'}:
     \  x\in X,\, \|x\|\le 1\right\}  
$$
It is not difficult to see that the space
${N}\!{^p}\,(X,Y)$ is the image of
$ X^*\wh\ot^p Y$ under the natural mapping
$ X^*\wh\ot^p Y\to L(X, Y).$ For $ z\in  X^*\wh\ot^p Y,$
the induced operator we will denote by
$ \wt z.$

We shall need below the following reformulation of the corollary
1.2 from [2]:
\small

 $(*)$ There exists a reflexive separable Banach space
$ E$ such that for each $ r\neq 2$ the canonical mapping
$ E^*\wh\ot^r E\to N^r(E, E)$ is not one--to--one.
  \small

The dual space of $X^*\wh\ot^p Y $ is isometrically equal to
the space $ \Pi_{p'}^d(Y, X^{**}),$ where by
$ \Pi_{p'}^d$ we denote the ideal, dual to the ideal of
absolutely
$ p'$-summing operators with the corresponding norm (see [1]):
if $z\in X^*\wh\ot^p Y $ and  $U\in \Pi_{p'}^d(Y, X^{**}),$
then the duality defined with the help of trace:\,
$ \tr U\circ z.$ Note that in case when one of the space
$ X^*$ or $ Y$ possesses the AP, the canonical mapping
$ X^*\wh\ot^p Y\to L(X, Y)$ is one--to--one and thus
we can write in this case:
  $ X^*\wh\ot^p Y= {N}\!{^p}\,(X,Y).$

The assertion $ (*)$ can be reformulated by the following way:
\small

$(**)$ There exists a reflexive separable Banach space
$ E$ such that for each $ r\neq 2$ one can find
a tensor $ z\in E^*\wh\ot^r E$ and an operator
$ U\in \Pi^d_{r'}(E,E)$ for which
 $ \tr U\circ z=1$ and the associated with
$ z$ operator $ \wt z=0.$
      \small

\proclaim {Theorem {\rm 1}}
Let $\, p\in [1,+\infty],$ $ T\in L(X,Y)$
and either
$\, X^*\in \,AP\ $ or $\, Y^{***}\in \,AP.$
If $ T\in N^p(X, Y^{**}),$
then $T\in N^p(X,Y).$ In other words, under these conditions
from the $ p$-nuclearity of the conjugate operator $ T^*$ it follows
that the operator $ T$ belongs to the space $ N^p(X, Y).$
\endproclaim

\demo{Proof}
Suppose there exists such an operator
    $ T\in L(X,Y),$
that
$ T\notin N^p(X,Y),$ but $ \pi_Y\,T\in N^p(X,Y^{**}).$
Since either $ X^*$ or $ Y^{**}$  has the $ AP,$
 $ N^p(X,Y^{**})=X^*\widehat\otimes^p Y^{**}.$
Therefore the operator $ \pi_Y\,T$ can be identified with the
tensor element
$ t\in X^*\widehat\otimes^p Y^{**};$
in addition, by the choice of $ T,$ \
$  t\notin X^*\widehat\otimes^p Y$ \ (the space
$  X^*\widehat\otimes^p Y$
is considered as a subspace of the space
$  X^*\widehat\otimes^p Y^{**}$\!).
Hence there is an operator
$ U\in \Pi^d_{p'}(Y^{**},X^{**})=\( X^*\widehat\otimes^p Y^{**}\)^*,$
with the property that
$ \tr U\circ t=\tr \(t^*\circ \( U^*|_{X^*}\) \)=1$ and
$ \tr U\circ \pi_Y\circ z=0$ for each
$ z\in X^*\widehat\otimes^p Y.$
>From the last it is follows that, in particular, $ U\pi_Y=0$ and
$ \pi_Y^*\,U^*|_{X^*}=0.$
In fact, if
$ x'\in X^*$ and $ y\in Y,$ then
$$ <U\pi_Y\,y,x'> = <y, \pi_Y^*\,U^*|_{X^*}x'>
       = \tr \,U\circ (x'\otimes \pi_Y(y))=0.
$$
Evidently, the tensor element
$ U\circ t$  induces the operator
$ U\pi_Y T,$ which is equal identically to zero.

If $ X^*\in AP$ then
 $ X^*\widehat\otimes^p Y^{**}= N^p(X,Y^{**})$
and, therefore, this tensor element is zero what is contradicted
to the equality
$ \tr\, U\circ t=1.$

Let now $ Y^{***}\in AP.$ In this case
$$ V:= \( U^*|_{X^*}\)\circ T^*\circ \pi_Y^*: \
       Y^{***}\to Y^* \to X^*\to Y^{***}
$$
uniquely determines a tensor element
$ t_0$ from the projective tensor product
$ Y^{****}\widehat{\otimes} Y^{***}.$
Let us take any representation $ t=\sum x'_n\otimes y''_n$ for $ t$
as an element of the space $ X^*\widehat{\otimes}^p Y^{**}.$
Denoting for the brevity the operator $ U^*|_{X^*}$ by $ U_*,$
we obtain:
$$\multline
   Vy'''=U_*\, \( T^*\pi_Y^*\,y'''\) =
    U_*\, \( (T^*\pi_Y^*\pi_{Y^*})\,\pi_Y^*\, y'''\) =
    U_*\, \( (\pi_Y T)^*\,\pi_{Y^*})\,\pi_Y^*\, y'''\) = \\ =
    U_*\, \( (\sum y''_n\otimes x'_n)\,\pi_{Y^*})\,\pi_Y^*\, y'''\)
    =U_*\, \( \sum <y''_n, \pi_Y^*\, y'''> \,x'_n \) = \\ =
      \sum <\pi_Y^{**}y''_n,  y'''> \,U_* x'_n.
  \endmultline
$$

So, the operator $ V$ (or the element $ t_0$) has in the space
$ Y^{****}\widehat\otimes Y^{***}$ the representation
$$ V= \sum \pi_Y^{**}(y''_n)\otimes U_* (x'_n).
$$
Therefore,
$$  \tr t_0=\tr V= \sum <\pi_Y^{**}(y''_n), U_* (x'_n)> =
        \sum <y''_n, \pi_Y^*\,U_* x'_n> =  \sum 0=0.
$$
On the other hand,
$$  Vy'''= U_* \( \pi_Y T\)^* y'''= U_*\circ t^* (y''')=
     U_*\, \( \sum <y''_n, y'''> \, x'_n\)=
      \sum <y''_n, y'''> \, U_* x'_n,
$$
whence $ V=\sum y''_n\otimes U_*(x'_n).$   Therefore
$$ \tr t_0=\tr V= \sum <y''_n, U_* x'_n> = \sum <Uy''_n, x'_n>
= \tr U\circ t=1.$$
The obtained contradiction completes the proof of the theorem.
 $\quad\blacksquare$

\enddemo

\proclaim {\bf Theorem 2}\it
For each $ r\in [1,\infty], r\neq 2,$ there exist a separable
space $ W$ and an operator $ T\in L(W^{**}, W)$ such that
$ W^{**}$ has a basis, $ T\in N^r(W^{**}, W^{**}),$ but
$ T\notin N^r(W^{**}, W).$
\endproclaim\rm

\demo{\it Proof}
Let us fix $ r\in [1,\infty], r\neq 2$ \, and take the triple $ (E, z,U)$
from the assertion $ (**).$
Let $ W$ be a separable space such that $ W^{**}$ has a basis
and there exists a linear homomorphism $ \ffi$ from $ W^{**}$
to $ E$ with the kernel
$ W\sbs W^{**}$ so that the subspace $ \ffi^*(E^*)$ in complemented
in $ W^{***}$ (see [4]). Lift the tensor element
$ z,$ lying in $ E^*\wh\ot^r E,$ up to an element \footnote{
If $ z=\sum_{k=1}^\infty \,e'_k\ot\,e_k $ is any representation of $ z$
in  $ E^*\wh\ot^r E,$  then we take $ \{ w''_n\}\sbs W^{**}$
in such a way that
the last sequence is absolutely $ r$-summing and
$ \ffi(w''_n)=e_n$ for every $ n.$}
$ \al\in E^*\wh\ot^r W^{**},$
so that $ \ffi\circ \al=z,$ and set $ V:= U\circ \ffi.$
Since $ \tr V\circ\al=\tr U\circ z=1$ and $ W^{**}$ has the
AP, then $ \wt \al=\al\neq 0.$
Besides, the operator $ \wt{\ffi\circ\al}:E\to W^{**}\to E,$
associated with the tensor $ \ffi\circ\al,$ is equal to zero. Therefore
$ \al(E)\subset \operatorname{ Ker}\ffi= W\sbs W^{**},$
that is the operator
$ \al$ is acted from $ E$ into $ W.$

Since the subspace $ \ffi^*(E^*)$ is complemented in $ W^{***},$
then $ \al\circ\ffi\in W^{***}\wh\ot^r W= N^r(W^{**},W)$ iff
$ \al\in E^{*}\wh\ot^r W= N^r(E,W).$

If $ \al\in  N^r(E,W),$
then for its arbitrary (nonzero!) $N^p$-representation
of the form $ \al=\sum e'_n\ot w_n$ the composition $ \ffi\circ\al$ is
a zero tensor element in
$ E^*\wh\ot E;$ but this composition
represents the element $ z, $ which, by its choice, can not be zero one.
Thus, $ \al\notin  N^r(E,W)$ and, thereby,
 $ \al\circ\ffi\notin W^{***}\wh\ot^r W= N^r(W^{**},W).$
On the other hand, certainly,
 $ \al\circ\ffi\in W^{***}\wh\ot^r W^{**}= N^r(W^{**},W^{**}).
 \quad\blacksquare$
\enddemo

So, the approximation conditions imposed on $X\,$ and $Y\,$
in the theorem 1 are essential.

It is clear that from $(**)$ and from the proof of Theorem 2
it follows that there exists a Banach space
$ W,$ for which the conclusion of the theorem is valid
for all mentioned values of the parameter $ r.$

\bigpagebreak

\centerline{REFERENCES}
\smallpagebreak

\ref \no 1 \by A. Pietsch \pages 536 p
 \paper  Operator ideals
 \yr 1978\vol
 \jour North-Holland
 \endref

\ref \no 2  \by O. I. Reinov \pages   125-134
\paper   Approximation properties of order p and the existence of
  non-p-nuclear operators with p-nuclear second adjoints
\yr 1982 \vol  109
\jour    Math. Nachr.
\endref

\ref \no 3  \by O. I. Reinov  \pages 277-291\nofrills
\paper
Approximation properties $ \operatorname{AP_s}$ and $p$-nuclear
operators {\rm(}the case when $ 0<s\le1)$
\yr 2000 \vol 270
\jour Zapiski nauchn. sem. POMI
\finalinfo (in Russian)
\endref

\ref \no 4\by J. Lindenstrauss \pages  279-284
\paper  On James' paper "Separable Conjugate Spaces"
\yr 1971\vol 9
\jour Israel J. Math.
\endref

\enddocument